\input amstex
\input epsf
\documentstyle{amsppt}
	\magnification=1200
	\rightheadtext{Circle packings in the unit disc}
	\hcorrection{0.25in}

\define\cp{circle packing}
\define\blp{Bl-packing}
\define\dbp{discrete Blaschke product}

\define\compk{\Bbb K}
\define\mcompk{$\compk$}
\define\vertk{{\Bbb K}^0}
\define\mvertk{${\Bbb K}^0$}

\define\medgek{${\Bbb K}^1$}
\define\facek{{\Bbb K}^2}
\define\mfacek{${\Bbb K}^2$}
\define\intk{\text{\rm{int}}\, {\Bbb K}^0}
\define\mintk{$\intk$}
\define\bdk{\text{\rm{bd}}\, {\Bbb K}^0}
\define\mbdk{$\bdk$}

\define\degc#1{\text{\rm{deg}}(\Bbb{#1})}
\define\mdegc#1{$\degc#1$}
\define\degv#1{\text{\rm{deg}}(#1)}
\define\mdegv#1{$\degv#1$}


\define\vcompk#1{{\compk}_{#1}}
\define\mvcompk#1{$\vcompk #1$}
\define\vsimpk#1#2{\vcompk #1^{#2}}
\define\mvsimpk#1#2{$\vsimpk #1#2$}

\define\vcomp#1#2{{\Bbb{#1}}_{#2}}
\define\mvcomp#1#2{$\vcomp #1#2$}


\define\interv#1{\text{\rm{int}}\, \vcompk #1^{0}}
\define\minterv#1{$\interv #1$}
\define\bdv#1{\text{\rm{bd}}\, \vcompk #1^{0}}
\define\mbdv#1{$\bdv #1$}

\define\p#1{\Cal{#1}}
\redefine\mp#1{$\p #1$}
\define\vp#1#2{\p #1_{#2}}
\define\mvp#1#2{$\vp #1#2$}

\define\smap#1{S_{\p #1}}
\define\msmap#1{$\smap #1$}

\define\carr#1{\text{\rm{carr}} (\p #1)}
\define\mcarr#1{$\carr #1$}


\define\brv#1{\text{\rm{br}}_V(\p #1)}
\define\mbrv#1{$\brv #1$}
\define\br#1{\text{\rm{br}}({\p #1})}
\define\mbr#1{$\br #1$}

\define\ord#1#2{\text{\rm{ord}}_{\Cal #1}(#2)}
\define\mord#1#2{$\ord #1#2$}

\define\fval#1{\text{\rm{val}}(#1)}
\define\mfval#1{$\fval #1$}
\define\pval#1{\text{\rm{val}}(\p #1)}
\define\mpval#1{$\pval #1$}

\def\today{\number\day\space\ifcase\month\or January \or
February \or March \or April \or May \or June \or July \or
August \or September \or October \or November \or
December\fi\space\number\year}


\topmatter
\title 
Circle packings in the unit disc
\endtitle
\author
Tomasz Dubejko
\endauthor
\address
Mathematical Sciences Research Institute, Berkeley, CA 94720\endaddress
\email
tdubejko$\@$msri.org\endemail
\thanks
Research at MSRI is supported in part by NSF grant DMS-9022140.\endthanks
\abstract
A Bl-packing is a (branched) circle packing that ``properly covers'' the unit disc.
We establish some fundamental properties of such packings.
We give necessary and sufficient conditions for their existence, prove their uniqueness, and show that their underlying surfaces, known as carriers, are quasiconformally equivalent to surfaces of classical Blaschke products.
We also extend the approximation results of [D1] to general combinatorial patterns of tangencies in Bl-packings.
Finally, a branched version of the Discrete Uniformization Theorem of [BSt1] is given.
\endabstract
\endtopmatter

\document
	\baselineskip=12pt
	\parindent=10pt
	\parskip=0pt

\heading Introduction \endheading

This is a continuation of the work initiated in [D1].
We are interested in (infinite) circle packings of the type considered in [D1], which we call Bl-packings.
These are (branched) \cp s that are contained in the closure of the unit disc $\bold D$, their circles do not accumulate in $\bold D$, and their boundary circles, if any, are internally tangent to $\partial \bold D$.
The existence of finite Bl-packings was proved in [D1]; the existence of infinite ones will be proved in this paper.

Our overall aim is to describe such packings, show similarities that they share with branched surfaces of classical Blaschke products, and establish their distortion properties.
We first prove that large collections of finite Bl-packings possess, under some rather natural restrictions, uniform bounds on their functions of local distortion which are defined in terms of ratios of radii of neighboring circles (Lemma~3.1).
This result turns out to be very useful tool for our goals.
It implies that large families of finite discrete Blaschke products of [D1] are K-quasiregular (Corollary~3.2).
In particular, this allows to extend approximation schemes of [D1] to general combinatorics and to show that every classical Blaschke product can be approximated uniformly on compact subsets of $\bold D$ by discrete Blaschke products in combinatorially less restrictive setup then that of [D1] (Theorem~3.3).

Lemma~3.1 is also used in proving the existence of infinite Bl-packings (Theorem~4.1).
Combining this with results of [DSt] and [D3] we obtain a branched version of the Discrete Uniformization Theorem of [BSt1], which shows the consistency in circle packing theory and its analogy to the classical theory of Riemann surfaces. 

We finish our investigation of properties of Bl-packings with Theorem~4.4 which generalizes Lemma~3.1 and quasiregularity result for finite discrete Blaschke products to infinite Bl-packings and infinite discrete Blaschke products.
The quantitative statement in Theorem~4.4 contains the inverse result to the Discrete Schwarz Lemma of [BSt2].

\heading 1.~Preliminaries \endheading

All simplicial 2-complexes considered in this paper are assumed to be oriented and simplicially equivalent to triangulations of a disc in the complex plane $\bold C$.
They could be infinite, possibly with boundary.
If \mcompk\ is such a complex then \mvertk, \mintk, \mbdk, \medgek, and \mfacek denote the sets of vertices, interior vertices, boundary vertices, edges, and faces of \mcompk, respectively.
If $u$ and $v$ are neighboring vertices in a complex, shortly $u\sim v$, then $uv$ stands for the edge between them.
The {\it degree} of a complex \mcompk\ is defined by $\degc K := \sup_{v\in \vertk}\{ \degv v\}$, where \mdegv v is the number of vertices adjacent to $v$.
A set $\{(v_1,k_1),\dots,(v_m,k_m)\}$ of pairs of interior vertices of \mcompk\ and positive integers is called a {\it branch structure} for \mcompk\ if for every simple closed edge path $\Gamma$ in \mcompk\ the number of its edges is at least $3+\sum_{i=1}^m\text{{\rm ind}}_{\Gamma}(v_i)k_i$, where $\text{{\rm ind}}_{\Gamma}(v_i)$ is 1 when $v_i$ is enclosed by $\Gamma$ and 0 otherwise.

A collection $\{C(v)\}_{v\in \vertk}$ of circles in $\bold C$ is said to be a {\it circle packing} for a complex \mcompk\ if $\langle C(u),C(v),C(w) \rangle$ is a positively oriented triple of mutually and externally tangent circles in $\bold C$ whenever $\langle u,v,w \rangle$ is an oriented faces of \mcompk. 
A \cp\ $\p P = \{C_{\p P}(v)\}$ for \mcompk\ determines the simplicial map $\smap P:\compk \to \bold C$ by the condition that $\smap P(v)$ is equal to the center of $C_{\p P}(v)$; the geometric simplicial complex $\smap P(\compk)$ is called the {\it carrier} of \mp P, \mcarr P.
For each $\triangle = \langle u,v,w \rangle \in \facek$ let $\alpha_{\p P}(v, \triangle)$ denote the angle sum at $\smap P(v)$ in the Euclidean triangle $\smap P(\triangle)$.
If $v\in \intk$ then $\Theta_{\p P}(v):= \sum_{\triangle \in \facek} \alpha_{\p P}(v,\triangle)$ turns out to be a positive integer multiple of $2\pi$ and is called the ({\it interior}) {\it angle sum} at $v$ induced by \mp P; if $v\in \bdk$ then $\gamma_{\p P}(v):= \sum_{\triangle \in \facek} \alpha_{\p P}(v,\triangle)$ is called the {\it boundary angle sum} at $v$ induced by \mp P.

We define the {\it set of branch vertices} of \mp P as $\brv P:= \{ v\in \intk: \Theta_{\p P}(v) > 2\pi \}$ and the {\it branch set} of \mp P as $\br P:= \{(v,\ord Pv)\}_{v\in \brv P}$, where $\ord Pv :=\tfrac{\Theta_{\p P}(v)}{2\pi}-1$.
A \cp\ is said to be {\it univalent} if all its circles have mutually disjoint interiors, in which case its set of branch vertices must be empty.

Suppose \mp P and \mp Q are \cp s for a complex \mcompk.
Then the simplicial map $f:\carr P \to \carr Q$, $f(\smap P(v)) = \smap Q(v)$ for $v\in \vertk$, is called the {\it cp-map} (read, \cp\ map) from \mp P to \mp Q; in this case \mp P and \mp Q are called the {\it domain} and {\it range} packings of $f$, respectively.
The associated affine map $f^\# : \carr P \to (0,\infty)$ determined by $f^\#(\smap P(v)) = \tfrac{r_{\p Q}(v)}{r_{\p P}(v)}$, $v\in \vertk$, is called the {\it ration map} from \mp P to \mp Q, where $r_{\ast}(v)$ is the radius of the circle $C_{\ast}(v)$.
Notice that if \mp P is univalent then $f$ and $f^\#$ are functions defined in a domain of $\bold C$ (=\mcarr P).

Recall that the valence \mfval f of a function $f$ is the least upper bound on the number of elements in the preimage under $f$ of any given point from the range.
The {\it valence of a packing} $\p P=\{C_{\p P}(v)\}$ is defined by $\pval P := \sup \{ n: \bigcap_{i=1}^n D_{\p P}(v_i) \ne \emptyset \}$, where $D_{\p P}(v)$ is the closed disc bounded by $C_{\p P}(v)$.
It is easy to see that $\fval f\le \pval P$ for every cp-map $f$ whose domain packing is univalent and whose range packing is \mp P.

Finally, we will say that a map $f: \Omega \to \Omega$, $\Omega \subset \bold C$, is an automorphism of $\Omega$ if it is 1-to-1, onto, and conformal.
We note that all automorphisms of $\bold C$ are similarities of $\bold C$, and all automorphisms of $\bold D$ are M\"obius transformations fixing $\bold D$.

\heading 2~Bl-packings and discrete Blaschke products \endheading

We begin with the following definition.

\proclaim{Definition~2.1}
Let \mcompk\ be a (finite or infinite) triangulation of a disc.
A \cp\ \mp P for \mcompk\ is said to be a Bl-packing if \mp P is contained in $\overline{\bold D}$, its circles have no accumulation point in $\bold D$, and its boundary circles, if any, are internally tangent to $\partial \bold D$.
\endproclaim

\midinsert
\vskip4.9truecm{ \epsfysize=7.5truecm \epsffile[72 161 540 630]{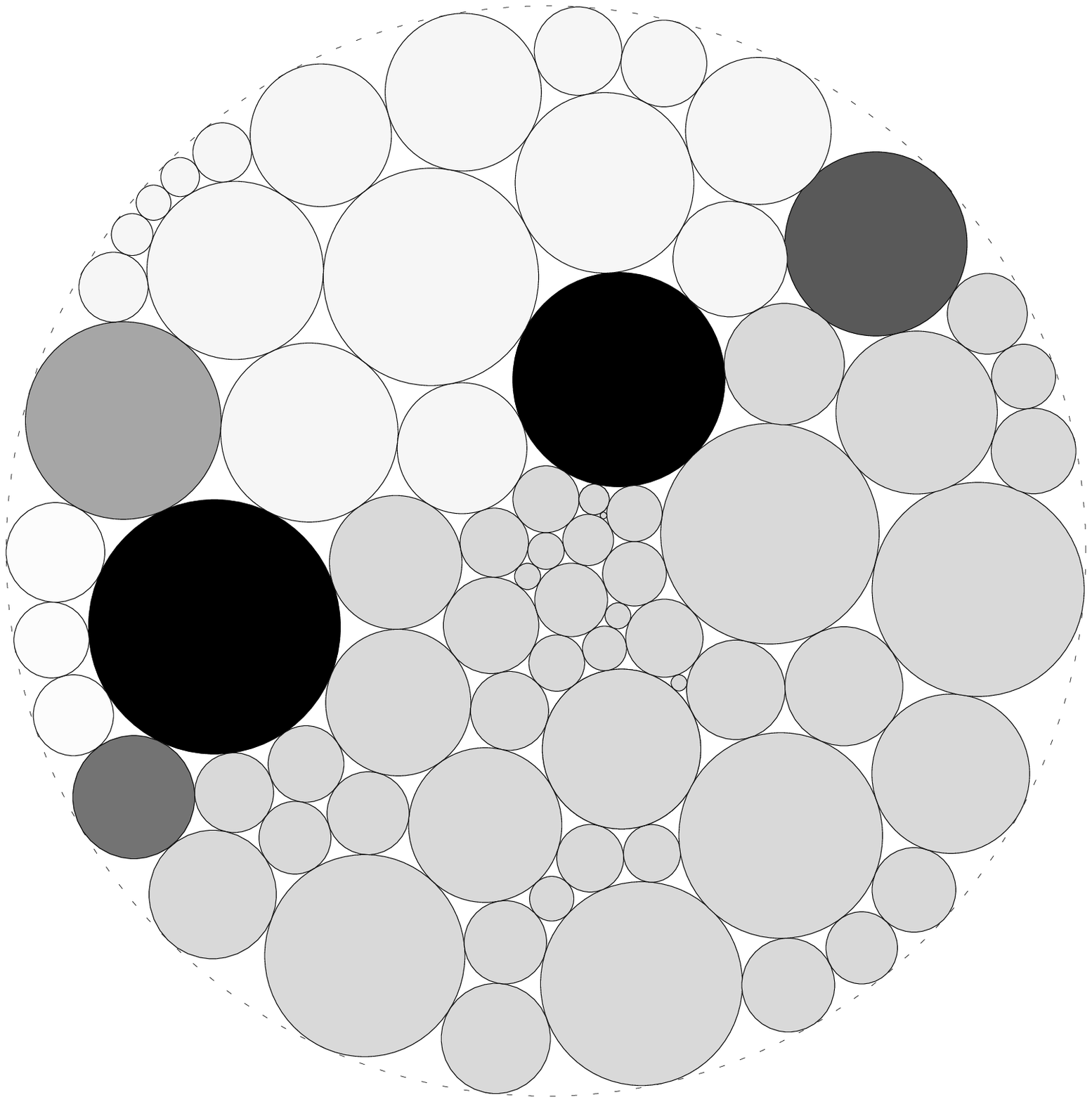}}
\vskip-1.0truecm\hskip0.25truecm (a)

\vskip-12.3truecm{ \hskip7.2truecm{\epsfxsize=7.5truecm \epsfysize=2.5truecm 
\epsffile[72 161 540 630]{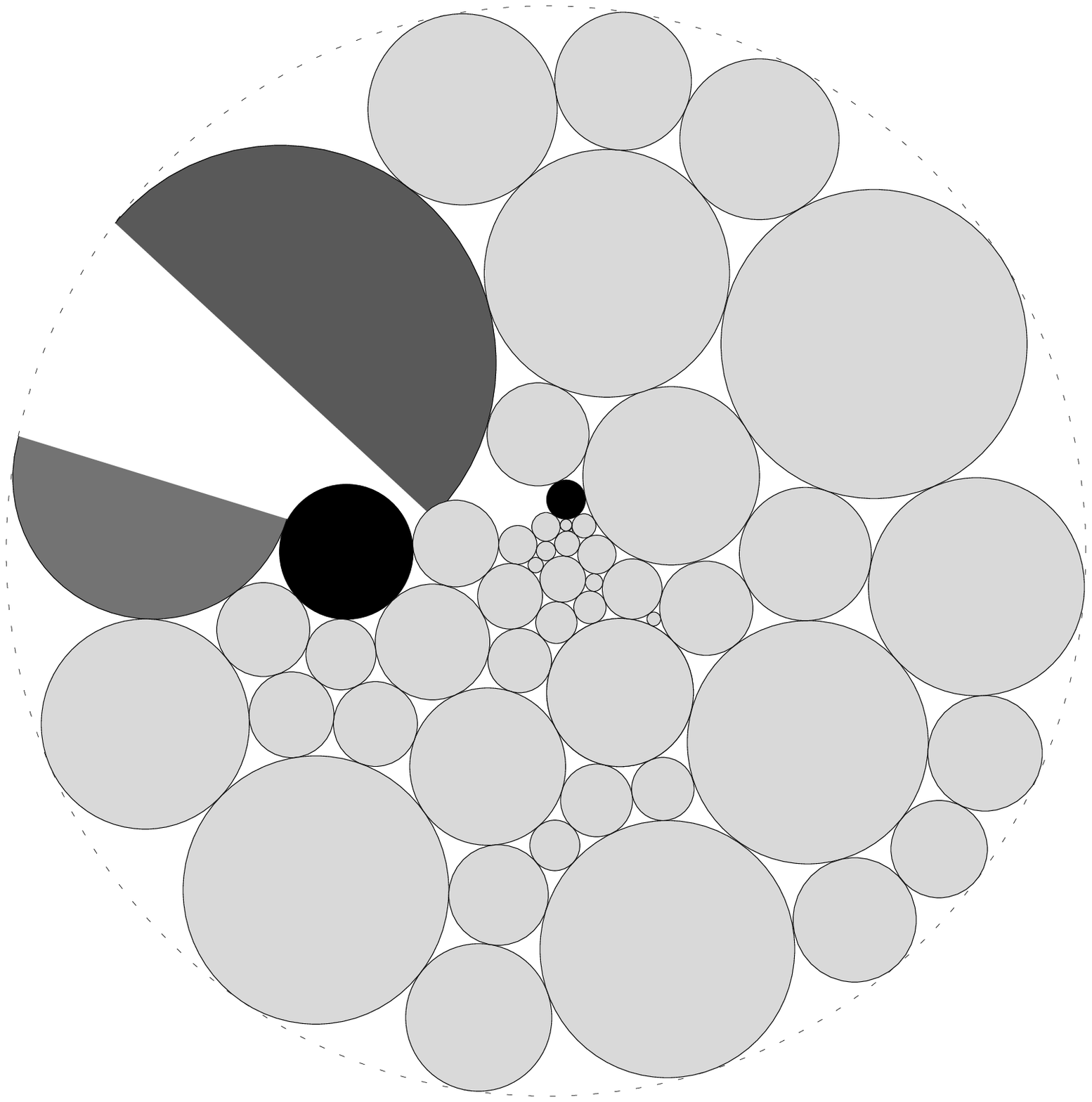}}}

\vskip0.02truecm{ \hskip7.2truecm{\epsfxsize=7.5truecm \epsfysize=2.5truecm \epsffile[72 161 540 630]{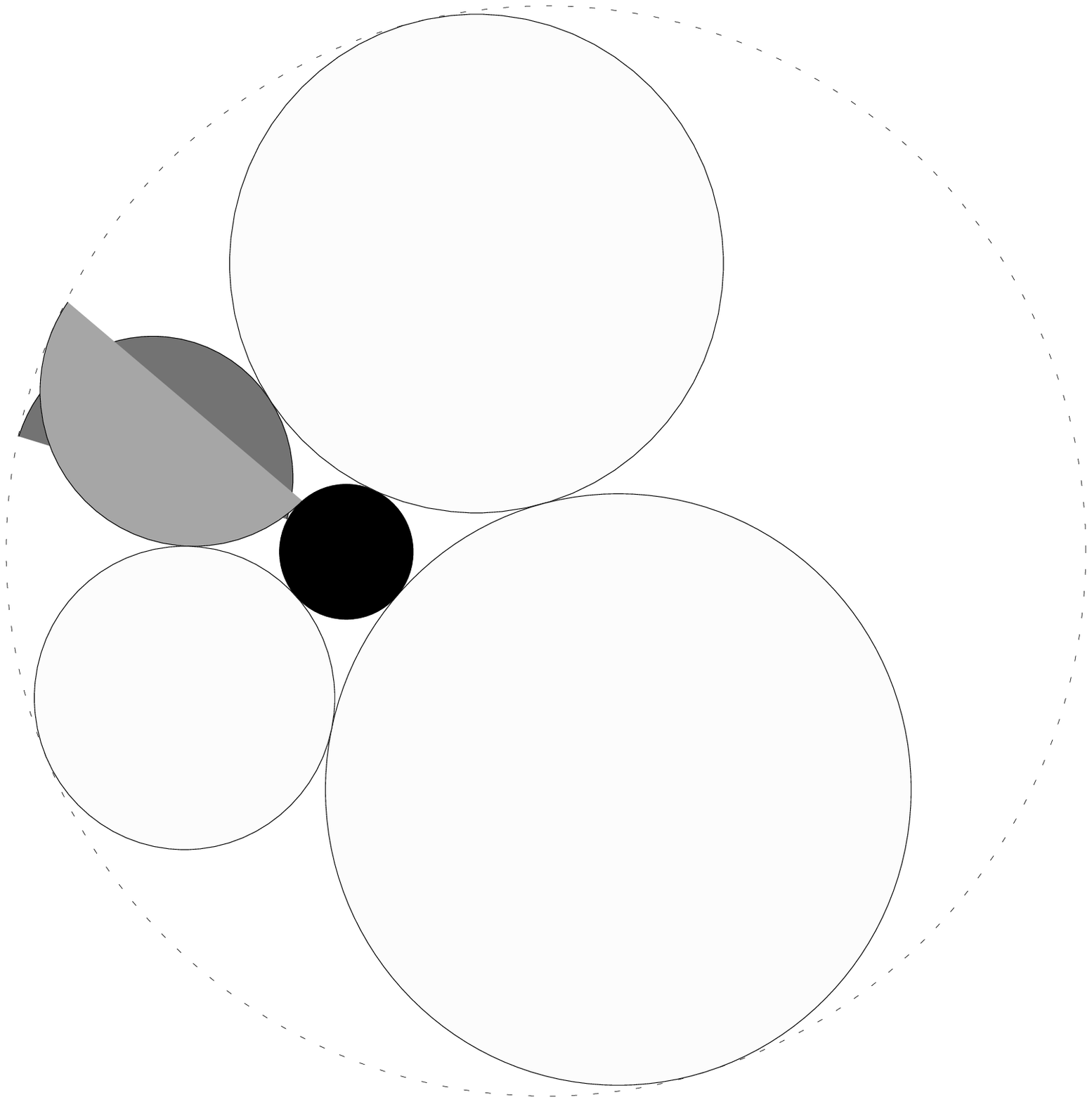}}}

\vskip0.02truecm{ \hskip7.2truecm{\epsfxsize=7.5truecm \epsfysize=2.5truecm \epsffile[72 161 540 630]{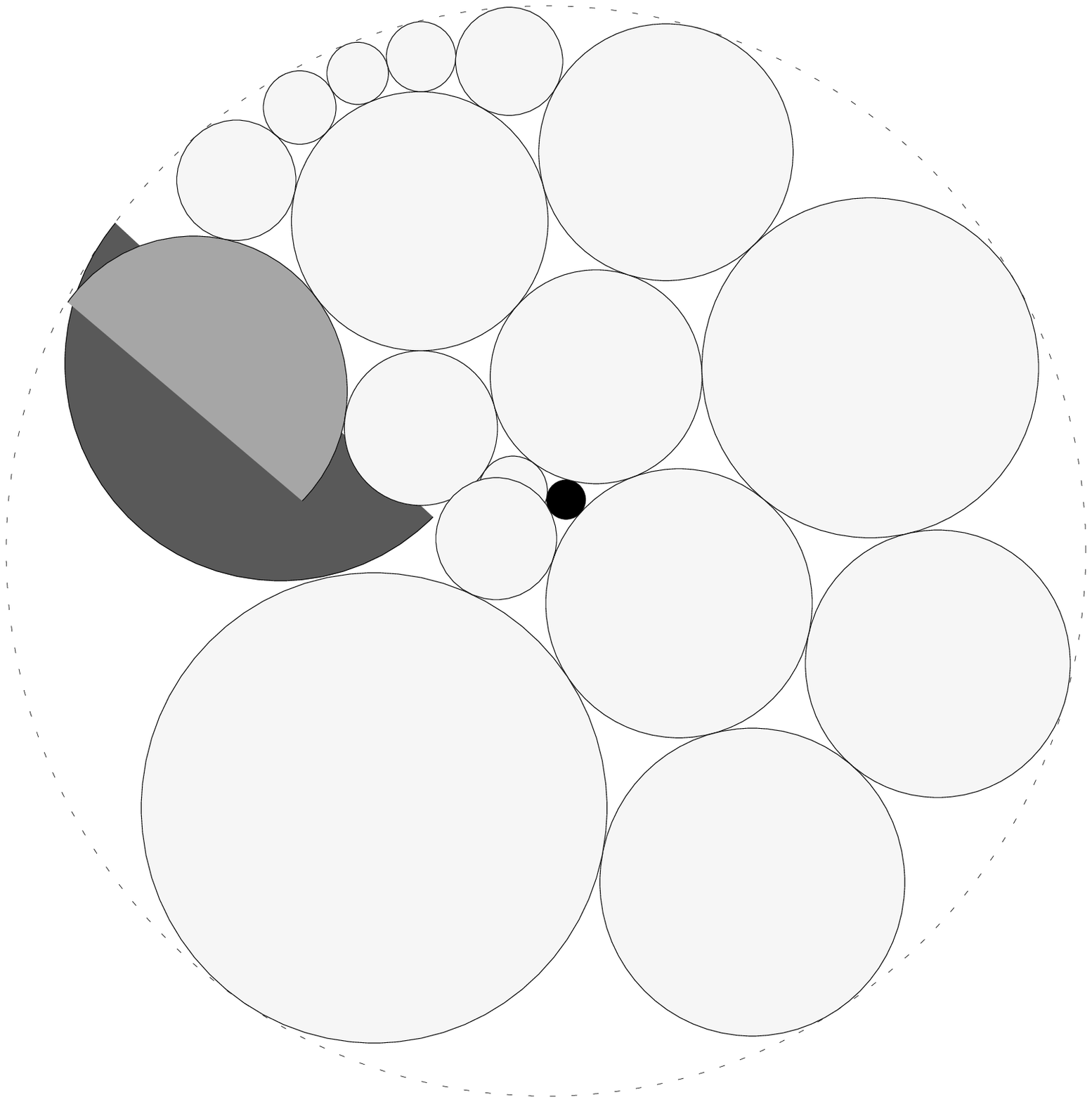}}}

\vskip0.02truecm\hskip9.4truecm\vbox{\epsfxsize=3.0truecm \epsfysize0.5truecm\epsffile[263 349 346 445]{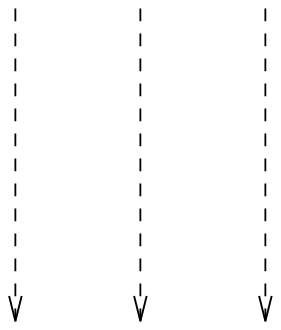}}

\vskip0.02truecm\hskip7.2truecm{\epsfysize=7.5truecm \epsffile[72 161 540 630]{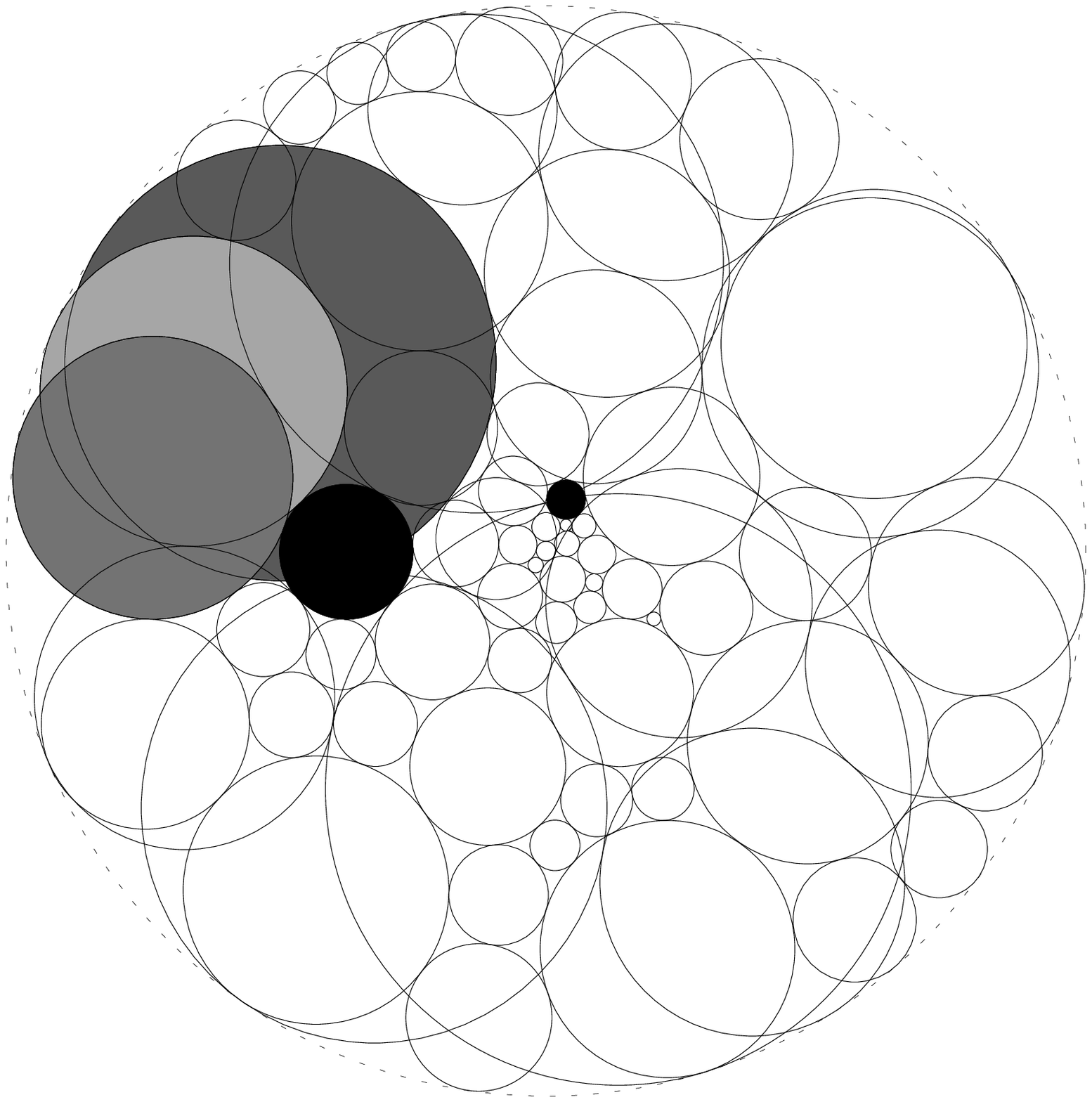}}
\vskip-1.0truecm\hskip13.5truecm (b)


\captionwidth{25pc}
\botcaption{Figure~1\footnotemark"*"} Bl-packings for the same complex: (a)~univalent and (b)~3-fold branched (with black circles corresponding to branch vertices of order 1) and its decomposition into locally univalent sheets.
\endcaption
\endinsert
\footnotetext"*"{This figure was created with the help of {\tt CirclePack} by the kind permission of Ken Stephenson}%

The existence of finite Bl-packings (i.e., Bl-packings for finite complexes) was proved in [D1].
The existence of infinite ones is more subtle and will be shown in Section~4.
Leaving aside the issue of existence we now establish some fundamental properties of Bl-packings.

\proclaim{Proposition~2.2} 
Suppose \mp P is a Bl-packing for \mcompk.
Then the following hold.
\roster
\item
There exists a univalent \blp\ for \mcompk.
\item
If $\tilde \p P$ is a univalent Bl-packing for \mcompk\ then there exists an extension $\bar f:\bold D \to \bold D$ of the cp-map $f$ from $\tilde \p P$ to \mp P such that $\bar f = \phi \circ h$, where $\phi: \bold D \to \bold D$ is a classical, finite Blaschke product and $h: \bold D \to \bold D$ is a homeomorphism.
Moreover, 
$$
\pval P = \text{\rm{val}}(\bar f) = 1 + \sum_{v\in \brv P} \ord Pv,
$$
and the branch set of $\phi$ is equal to $\{ (h(S_{\tilde \p P}(v)), \ord Pv) \}_{v\in \brv P}$.
\item
If \mp Q is \blp\ such that $\br Q = \br P$ then there exists a M\"obius transformation $M$ preserving $\bold D$ such that $\p Q = M(\p P)$, i.e. \blp s are uniquely determined, up to automorphisms $\bold D$, by their branch sets.
\endroster
\endproclaim

\demo{Proof}
(1)~This part follows easily, by switching to hyperbolic metric in $\bold D$, from the Perron-type arguments of [BSt2] (see [BSt2, Thm.~4]).

(2)~We first define a map $\bar f$.
If \mcompk\ has no boundary then we set $\bar f$ to be $f$. 
For \mcompk\ with boundary the construction of $\bar f$ is as follows.

Suppose $u,v\in \bdk$, $u\sim v$.
Define $\lambda_{\ast}(v) := C_{\ast}(v) \cap \partial \bold D$ and $\mu_{\ast}(u,v) := C_{\ast}(u) \cap C_{\ast}(v)$, where $\ast$ denotes \mp P or $\tilde \p P$.
Also, if there exists a component of $\bold D \setminus (C_{\ast}(u) \cup C_{\ast}(v))$ that has empty intersection with $\bigcup_{\triangle \in \facek : uv\subset \triangle} S_{\ast}(\triangle)$, we denote it by $A_{\ast}(u,v)$ (there could be at most one such a component).
Let $M(u,v)$ be the M\"obius transformation mapping $\overline{A_{\tilde \p P}(u,v)}$ onto $\overline{A_{\p P}(u,v)}$ with points $\lambda_{\tilde \p P}(u)$, $\lambda_{\tilde \p P}(v)$, and $\mu_{\tilde \p P}(u,v)$ mapped to $\lambda_{\p P}(u)$, $\lambda_{\p P}(v)$, and $\mu_{\p P}(u,v)$, respectively (Fig.~2).
Thus $f$ and the maps $M(u,v)$, $u,v\in \bdk$, $u\sim v$, induce a well-defined and continuous function $\hat f : \text{carr}(\tilde \p P) \cup (\bigcup_{v\in \bdk}C_{\tilde \p P}(v)) \cup (\bigcup\Sb u,v\in \bdk\\
u\sim v\endSb A_{\tilde \p P}(u,v)) \to \bold D$.
To obtain $\bar f$ we extend each $\hat f|_{C_{\tilde \p P}(v)}$, $v\in \bdk$, radially to the interior of $D_{\tilde \p P}(v)$.

\midinsert
\vskip0.0truecm\hskip6.0truecm $C_{\tilde \p P}(v)$\hskip6.0truecm $C_{\p P}(u)$
\epsfysize=6.8truecm
\vskip-0.5truecm\centerline{\epsffile[10 262 600 530]{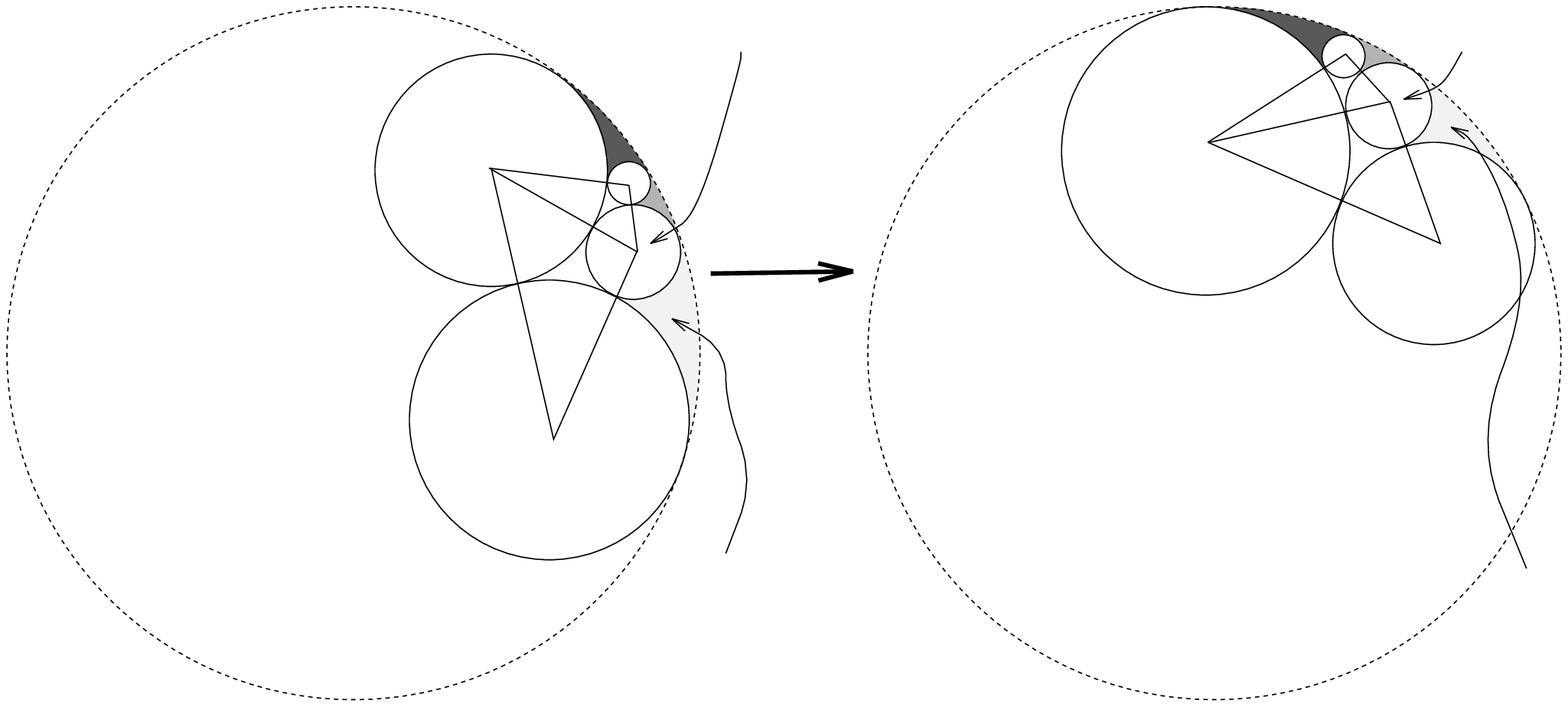}}
\vskip-4.9truecm\hskip7.0truecm $\bar f$
\vskip0.0truecm\hskip12.9truecm $C_{\p P}(u)$
\vskip1.5truecm\hskip4.5truecm $C_{\tilde \p P}(u)$
\vskip0.2truecm\hskip5.7truecm $A_{\tilde \p P}(u,v)$
\vskip-0.35truecm\hskip13.1truecm $A_{\p P}(u,v)$
\vskip0.0truecm\hskip2.8truecm $\tilde \p P$ \hskip8.0truecm $\p P$
\captionwidth{16pc}
\botcaption{Figure~2} The extension.
\endcaption
\endinsert

Now recall that a continuous function $g:X \to Y$, $X,Y\subset \bold C$, is proper if preimage of every compact set in $Y$ is compact.
Definitions of Bl-packings and cp-maps together with the above construction imply that $\bar f$ is a proper mapping.
Thus, Stoilow's theorem [LV] and the fact that an analytic proper function $\Phi:\bold D \to \bold D$ is a finite Blaschke product show that $\bar f = \phi \circ h : \bold D\to \bold D$, where $h : \bold D\to \bold D$ is a homeomorphism and $\phi : \bold D\to \bold D$ is a finite Blaschke product.
In particular, $\pval P = \text{\rm{val}}(\phi) = 1+\sum_{v\in \brv P} \ord Pv <\infty$.

(3)~The assertion of this part follows from the Discrete Schwarz and Distortion Lemmas of [DSt].
\qed
\enddemo

In [D1] the notion of finite, discrete Blaschke products was introduced.
In the light of the last proposition it is natural to extend this notion to the infinite case.

\proclaim{Definition~2.3}
A discrete Blaschke product is a cp-map between two Bl-packings with the domain packing being univalent.
\endproclaim

If $f$ is a \dbp\ then the extension of $f$ constructed in the proof of Proposition~2.3(3) will be denoted by $\bar f$.

\heading 3.~Finite \blp s and \dbp s: quasiregularity and approximation \endheading

In this section we are primarily interested in a local behavior of finite \blp s and \dbp s.
We establish here our key result, Lemma~3.1, which is followed by two applications.
The first application implies that \dbp s are quasiregular mappings.
The second one shows that approximation schemes of [D1] extend to general combinatorics. 

For the purpose of studying local behavior of \blp s we introduce a function that measures local distortion of packings.
If \mp P is a circle packing for \mcompk\ then the function $\mu_{\p P}: \vertk \to (0,\infty)$ defined by 
$$
\mu_{\p P}(v):= \max\Sb u\in \vertk \\ u\sim v\endSb \, \frac12 \left( \frac{r_{\p P}(u)}{r_{\p P}(v)} + \frac{r_{\p P}(v)}{r_{\p P}(u)} \right)
$$
will be called the {\it function of local distortion} of \mp P.

We can now state 

\proclaim{Lemma~3.1}
Let $d$ and $\ell$ be positive integers, and $\rho \in (0,1)$.
There exists a constant $\kappa = \kappa(d,\ell,\rho)\ge 1$, depending only on $d$, $\ell$, and $\rho$, such that if \mp P is a \blp\ for a complex \mcompk\ satisfying 1)~$\degc K \le d$, 2)~$\pval P \le \ell$, and 3)~$\smap P(\intk) \cap \{|z|\le \rho \} \ne \emptyset$ and $\{\smap P(v) \}_{v\in \brv P} \subset \{|z|\le \rho \}$ then 
$$
\mu_{\p P}(v)\le \kappa \quad \text{for all $v\in \vertk $}.
\tag{*}
$$
\endproclaim

The importance of this lemma lies in the fact that (*) holds for all vertices in the complex (cf. [D1, Lemma~6.1] and [D2, Corollary~4.9]).
Its proof will be given in the last section.
Now couple remarks regarding its hypotheses.
It is immediate that for univalent \blp s only the conditions~1) and the first part of 3) are necessary; the remaining conditions are satisfied anyway in this case.
For \blp s which are not univalent it is not hard to see that all conditions~1)~--~3) are necessary, though the first part of the condition~3) can be removed because it follows from the second part.

The following is almost an immediate consequence of Lemma~3.1.

\proclaim{Corollary~3.2}
Let $d$, $\ell$, and $\rho$ be as in Lemma~3.1.
Then there exists a constant $K= K(d,\ell,\rho) \ge 1$, depending only on $d$, $\ell$, and $\rho$, such that if $f$ is a \dbp\ whose domain and range packings  satisfy 1)--3) of Lemma~3.1 then the corresponding extension map $\bar f: \bold D \to \bold D$ is $K$-quasiregular.
\endproclaim
\demo{Proof}
From (*) it follows that there is a constant $\kappa' = \kappa'(d,\ell,\rho) \ge 1$ such that for every $f$ satisfying hypotheses of the corollary $\tfrac 1{\kappa'}\alpha_{\p R_f}(v,\triangle) \le \alpha_{\p D_f}(v,\triangle) \le \kappa'\alpha_{\p R_f}(v,\triangle)$ for all $v\in \vertk$ and $\triangle = \langle u,v,w \rangle \in \facek$, where $\p D_f$ and $\p R_f$ are the domain and range packings of $f$, respectively, and \mcompk\ is their complex.
Hence, the definition of \dbp s implies that there is $K_1 = K_1(d,\ell,\rho)$ such that $\bar f|_{\text{\rm{carr}} (\p D_f)} = f$ is $K_1$-quasiregular. 
It is also easy to see that (*) and the definition of extensions of \dbp s imply that there exists $K_2 = K_2(d,\ell,\rho)$ such that $\bar f|_{\bold D \setminus \text{\rm{carr}} (\p D_f)}$ is $K_2$--quasiregular.
By taking $K:= \max_{i=1,2}\{ K_i\}$ we obtain the assertion of the corollary.
\qed
\enddemo

In [D1] a constructive method using hexagonal complexes was developed to show that every classical Blaschke product can be approximated uniformly on compact subsets of $\bold D$ by \dbp s.
Corollary~3.2 allows for a generalization of this approximation scheme to other then hexagonal combinatorics.

\proclaim{Theorem~3.3}
Let $\phi$ be a finite Blaschke product.
Write $\{ (x_1,k_1),\dots ,(x_m,k_m) \}$ for the branch set of $\phi$.
Suppose $\{ \Cal D_n \}$ and $\{ \Cal R_n \}$ are sequences of finite \blp s such that
\roster
\item
$\Cal D_n$ and $\Cal R_n$ are packings for the same complex $\compk_n$, and $\Cal D_n$ is univalent,
\item
there exists $d$ such that $\text{\rm{deg}}(\compk_n) \le d$ for all $n$,
\item
$\lim_{n\to \infty} \sigma_n =0$, where $\sigma_n := \sup_{v\in \vertk_n}\{r_{\Cal D_n}(v)\}$,
\item
for each $n$ there are vertices $\dot v_n, \ddot v_n\in \vertk_n$ such that $S_{\Cal D_n}(\dot v_n)=0$, $S_{\Cal R_n}(\dot v_n) = \phi(0)$, $\lim_{n\to \infty}S_{\Cal D_n}(\ddot v_n) = 1/2$, and $\text{\rm{Arg}}(S_{\Cal R_n}(\ddot v_n)) = \text{\rm{Arg}}(\phi(1/2))$, where $\text{\rm{Arg}}(z)$ denotes the argument of the complex number $z$,
\item
each $\text{\rm{br}}_V(\Cal R_n)$ can be written as a disjoint union of sets $V_1(n),\dots , V_m(n)$ such that for each $n$ one has $\sum_{v\in V_i(n)} \text{\rm{ord}}_{\Cal R_n}(v) = k_i$, $i=1,\dots ,m$, and $\lim_{n\to \infty} S_{\Cal D_n}(V_n(i)) = \{x_i\}$, $i=1,\dots ,m$, as the set limit.
\endroster 
Denote by $f_n$ the \dbp\ from $\Cal D_n$ to $\Cal R_n$.
Then the functions $f_n$ and $f^{\#}_n$ converge uniformly on compact subsets of $\bold D$ to $\phi$ and $|\phi'|$, respectively.
\endproclaim

\demo{Proof}
The proof will be the same as the proof of Theorem~7.2 in [D1] except for two modifications.
The first modification is required in order to show that $\lim_{n\to \infty} \delta_n =0$, where $\delta_n := \sup_{v\in \vertk_n}\{r_{\Cal R_n}(v)\}$.
The second one is needed in proving that $\lim_{n\to \infty} f_n$ is a holomorphic function.

By applying if necessary the automorphism $\Phi$ of $\bold D$ such that $\Phi(\phi(0))=0$ and $\text{\rm{Arg}}(\Phi(\phi(1/2))) = \text{\rm{Arg}}(\phi(1/2))$ to packings $\p R_n$, we can assume that $\phi(0)=0$.

Recall [LV] that a $K$-quasiconformal homeomorphism $g:\bold D \to \bold D$, $g(0)=0$, satisfies
$$
|g(z_1)-g(z_2)|\le 16|z_1-z_2|^{1/K} \quad \text{for all $z_1,z_2\in \overline{\bold D}$}.
\tag{3.1}
$$
Recall also that for each $n$ the extension map $\bar f_n$ of $f_n$ has a decomposition $\bar f_n = \phi_n \circ h_n$, where $\phi_n$ is a Blaschke product and $h_n$ is a self-homeomorphism of $\bold D$, $h_n(0)=0$, $h_n(1/2)\in (0,1)$, and $\phi_n(0)=0$ (because $\phi(0)=0$).

The condition~(5) implies that each $\phi_n$ is a Blaschke product of degree $\ell = 1+ \sum_{i=1}^m k_i$; in particular $\text{\rm{val}}(\p R_n)\le \ell$ for every $n$.
Moreover, (5) and the Discrete Schwarz Lemma of [BSt2] imply that there is $\rho$ such that $S_{\Cal D_n}(\text{\rm{br}}_V(\p R_n))$ and $S_{\Cal R_n}(\text{\rm{br}}_V(\p R_n))$ are contained in $\{|z|\le \rho\}$ for all $n$.
Therefore, by Lemma~3.1, the functions $h_n$ are $K$-quasiconformal and form a normal family whose every convergent subsequence has a $K$-quasiconformal self-homeomorphism of $\bold D$ as its limit [LV, Ch.~II, Thm.~5.5].
This gives that $\{\phi_n\}$ is not only a normal family but it is also compact in the class of Blaschke products, and it follows from (3.1) and the equicontinuity of normal families that $\lim_{n\to \infty} \delta_n =0$.

Thus we are now in a position of copying all quasiregular arguments (except for the argument involving conformality) used in the proof of Theorem~7.2 in [D1] to our current setting.
By replacing the hexagonal complexes $H(n)$ of [D1] by the complexes $\compk_n$ and by repeating the arguments of [D1] we arrive to the following conclusion:
from every subsequence of $\{f_n\}$ we can choose further subsequence $\{f_{n_l}\}$ such that: 1)~$h_{n_l} \to h$ uniformly on compacta of $\bold D$ as $l\to \infty$, where $h:\bold D \to \bold D$ is a $K$-quasiconformal homeomorphism, $h(0)=0$, and $h(1/2)\in (0,1)$, and 2)~$\phi_{n_l} \to \psi$ uniformly on compacta of $\bold D$ as $l\to \infty$, where $\psi$ is a Blaschke product with branch set $\{(h(x_1),k_1),\dots , (h(x_m),k_m)\}$.

We will now show that $h$ is conformal; this involves the second modification.
Let $z\in \bold D \setminus \{x_1,\dots , x_m\}$.
The uniform convergence of $\{h_{n_l}\}$ and $\{\phi_{n_l}\}$ implies that there is $\epsilon>0$ and $L$ such that $f_{n_l}|_{B_{\epsilon}(z)}$ is 1-to-1 for every $l>L$, where $B_{\epsilon}(z)=\{\zeta : |\zeta -z|<\epsilon\}$.
In particular, the portion of $\p R_{h_l}$ corresponding to the portion of $\p D_{h_l}$ contained in $B_{\epsilon}(z)$ is univalent.
Since $\lim_{n\to \infty} \delta_n =0$, Theorem~2.2 of [HR] shows that the quasiconformal distortion of $f_{n_l}$ at $z$ goes to 0 as $l\to \infty$.
Hence we obtain that $h$ is 1-quasiconformal at $z$.
Thus $h$ is 1-quasiconformal in $\bold D \setminus \{x_1,\dots , x_m\}$, and it follows that $h$ is conformal in $\bold D$.

By Riemann mapping theorem, $h$ is the identity function and $\psi = \phi$.
Now standard arguments imply that the full sequence $\{f_n\}$ converges uniformly on compact subsets of $\bold D$ to $\phi$.

The convergence of $\{f^{\#}_n\}$ to $|\phi'|$ is a consequence of Theorem~1 of [DSt].
\qed
\enddemo

\heading 4.~Existence and properties of (infinite) Bl-packings \endheading

As we have mentioned earlier, the author gave in [D1] the necessary and sufficient conditions, Theorem~2.5, for the existence of finite \blp s. With the help of Lemma~3.1 we are now able to extend this result to the infinite case, and show the existence of infinite \blp s.

\proclaim{Theorem~4.1}
Let \mcompk\ be a triangulation of a disc with $\degc K < \infty$.
Then there exists a \blp\ for \mcompk\ with branch set $\goth B$ if and only if there exists a univalent \blp\ for \mcompk\ and $\goth B$ is a branch structure for \mcompk.
\endproclaim

\demo{Proof}
Since the necessity is a straightforward consequence of Proposition~2.2 and Theorems~2.2 of [D1], we only need to be concerned with the sufficient condition in our assertion and, due to [D1, Thm.~2.5], only when \mcompk\ is infinite.

Let $\compk_1\subset \compk_2\subset \dots$ be a sequence of subcomplexes of \mcompk\ such that $\{b_1,\dots ,b_m\}\subset \compk_1$, each $\compk_j$ is a finite triangulation of a disc, and $\bigcup_j\compk_j = \compk$.
Let $u_0$ and $u_1$ be two designated vertices in $\compk_1$.
Denote by $\p P_j$ the univalent \blp\ for $\compk_j$ such that $S_{\p P_j}(u_0)=0$ and $S_{\p P_j}(u_1)\in (0,1)$.
It follows from [BSt1, Lemma~5] and the hypothesis that the sequence of packings $\{\p P_j\}$ has the geometric limit which is a \blp\ $\p P$ for \mcompk\ such that $S_{\p P}(u_0)=0$ and $S_{\p P}(u_1)\in (0,1)$.
In particular, there is $\rho$, $0<\rho<1$, such that $S_{\p P_j}(\{b_1,\dots ,b_m\})\subset \{|z|<\rho\}$ for every $j$, and for each $v\in \vertk$, $\lim_{n\to \infty} r_{\p P_j}(v) = r_{\p P}(v)$.
By the hypothesis and [D1, Thm.~2.5], for each $j$ there exists a \blp\ $\ Q_j$ for $\compk_j$ such that $\text{\rm{br}}(\p Q_j) = \goth B$, $S_{\p Q_j}(u_0)=0$ and $S_{\p Q_j}(u_1)\in (0,1)$.
We will show that the sequence $\{\p Q_j\}$ has the geometric limit.

From Theorem~4 of BSt1 it follows that $S_{\p Q_j}(\{b_1,\dots ,b_m\})\subset \{|z|<\rho\}$ for every $j$.
Write $f_j$ for the \dbp\ from $\p P_j$ to $\p Q_j$.
Let $\bar f_j$ be the corresponding extension.
Then, by Corollary~3.2, $\{\bar f_j\}$ is a family of $K$-quasiregular mappings.
Moreover, for each $j$, $\bar f_j = \phi_j \circ h_j$, where $h_j:\bold D \to \bold D$ is $K$-quasiconformal homeomorphism, $h_j(S_{\p P_j}(u_0)) =h_j(0)=0$, $h_j(S_{\p P_j}(u_1))\in (0,1)$, and $\phi_j$ is a Blaschke product, $\phi_j(0)=0$.
Since $h_j$'s are uniformly bounded, $\{h_j\}$ is a normal family.
Thus, by taking a subsequence if necessary, we may assume that $\{h_j\}$ converges uniformly on compacta of $\bold D$ to a $K$-quasiconformal homeomorphism $h$ of $\bold D$ onto $\bold D$ (see [LV, Ch.~II, Thm.~5.5]).
In particular, there exists $\bar \rho$, $0<\bar \rho<1$, such that $h_j(S_{\p P_j}(\{b_1,\dots ,b_m\}))\subset \{|z|<\bar \rho\}$ for all $j$, i.e. the critical points of each $\phi_j$ are in $\{|z|<\bar \rho\}$.
Similarly, since $\{\phi_j\}$ is a normal family of Blaschke products of degree $1+\sum_{i=1}^m k_i$ and 0 is a fixed point for all $\phi_j$'s, we can assume that $\{\phi_j\}$ converges uniformly in $\overline{\bold D}$ to a Blaschke product $\phi :\overline{\bold D} \to \overline{\bold D}$ whose critical set is equal to $\bigl\{\bigl( h(S_{\p P}(b_j)), k_i\bigr)\bigr\}_{i=1}^m$.
Hence the sequence $\{\p Q_j\}$ has the geometric limit \mp Q which is a \cp\ for \mcompk\ such that $\phi \circ h$ is the cp-map from \mp P to \mp Q.
Thus \mp Q is the sought after \blp.
\qed
\enddemo

We will now combine the above result with a similar one proved in [D3] for \cp s in $\bold C$ to show the consistency in circle packing theory and its analogy to the classical theory of (branched) surfaces.
For this purpose we need to introduce a new term.

Recall that a cp-map $f$ is said to be a discrete complex polynomial (see [D2]) if the domain packing of $f$ is univalent and there exists a decomposition $f=\psi \circ h$, where $\psi$ is a complex polynomial and $h$ is self-homeomorphism of $\bold C$.
A \cp\ which is the range packing of a discrete complex polynomial will be called a {\it Pl-packing}.

By putting together Corollary~4.3 of [D3], Proposition~2.2, and Theorem~3.3 we obtain the following branched version of the Discrete Uniformization Theorem of [BSt1]:

\proclaim{Uniformization}
If \mcompk\ is a triangulation of a disc with $\degc K < \infty$ then one and only one of the following holds:
\roster
\item
for every branch structure $\goth B$ for \mcompk\ there exists a \blp\ for \mcompk\ with branch set $\goth B$;
\item
for every branch structure $\goth B$ for \mcompk\ there exists a  Pl-packing for \mcompk\ with branch set $\goth B$.
\endroster
Moreover, packings in (1) and (2) are determined uniquely up to automorphisms of $\bold D$ and $\bold C$, respectively.
\endproclaim

\remark{Remark~4.2}
The reader can notice that the above statement is not the complete branched analog of the Discrete Uniformization Theorem of [BSt1]; the case which we have left out is that of triangulations of the 2-sphere.
We do not know of any satisfactory answer to a question of the existence and uniqueness of branched \cp s for triangulations of the 2-sphere.
The best result so far is due to [BoSt] which gives the existence and uniqueness of such packings assuming that a half of the total branching occurs at one of vertices of a triangulation.
\endremark

\remark{Remark~4.3}
We expect the above result to be true for every triangulation of a disc, not necessarily with bounded degree.
In the univalent case, it was shown in [HSc] that if \mcompk\ is a triangulation of a disc ($\degc K \le \infty$) then there exits either a univalent Pl-packing for \mcompk\ or a univalent Bl-packing for \mcompk, but not both.
Moreover, in [D2] it was shown that if \mcompk\ has boundary and $\goth B$ is a branch structure for \mcompk\ then there exists a \cp\ for \mcompk\ contained in $\bold D$ with branch set $\goth B$, hence the Perron-type arguments of [BSt2] imply the existence of a \blp\ for \mcompk\ with branch set $\goth B$.
\endremark

We finish this section by extending the results from \S3 to infinite \blp s.

\proclaim{Theorem~4.4}
Suppose \mcompk\ is a finite or infinite triangulation of a disc, $\tilde \p P$ and \mp P are \blp s for \mcompk, and $\tilde \p P$ is univalent.
Let $d$ and $\ell$ be positive integers, and $\rho \in (0,1)$.
In addition to the constants $\kappa(d,\ell,\rho)$ and $K(d,\ell,\rho)$ of Lemma~3.1 and Corollary~3.2, respectively, there exist constants $\lambda =\lambda(d,\ell,\rho)$ and $\eta =\eta(d,\ell,\rho)$, depending only on $d$, $\ell$, and $\rho$, such that if $\degc K\le d$, $\pval P \le \ell$, and $\smap P(\brv P)\subset \{|z|\le \rho\}$ then:
\roster
\item
$\mu_{\p P}(v)\le \kappa$ for all $v\in \vertk$;
\item
if $S_{\tilde \p P}(\intk) \cap \{|z|\le \rho \} \ne \emptyset$ and $f$ is the \dbp\ from $\tilde \p P$ to \mp P then the corresponding extension map $\bar f:\bold D\to \bold D$ is $K$-quasiregular;
\item
if $S_{\tilde \p P}(v_0)=\smap P(v_0)=0\in \bold D$ for some $v_0\in \intk$ then
$$
\lambda (r_{\tilde \p P}(v_0))^{\eta}\le r_{\p P}(v_0)\le r_{\tilde \p P}(v_0).
$$
\endroster
\endproclaim

\demo{Proof}
If \mcompk\ is finite then (1) and (2) are just restatements of Lemma~3.1 and Corollary~3.2.
If \mcompk\ is infinite then the assertions in (1) and (2) follow Lemma~3.1, Corollary~3.2, and limit-type arguments similar to those used in the proof of Theorem~4.1.

We will now prove (3).
The inequality $r_{\p P}(v_0)\le r_{\tilde \p P}(v_0)$ is a consequence of the Discrete Schwarz Lemma of [DSt] (see also [BSt2, Thm.~4]).
To show the other inequality in (3) we refer to quasiregular arguments.
Recall that if $g:\Omega\to \bold C$, $\Omega\subset \bold C$, is a $K$-quasiregular mapping, $A\subset \Omega$ is open, bounded, $\overline A \subset \Omega$, and $\partial g(A) = g(\partial A)$, then for every  compact subset $B$ of $A$,
$$
\text{\rm{cap}}(A,B) \le 
K\, \fval g \, \text{\rm{cap}}(g(A),g(B))
\tag{4.1}
$$
where $\text{\rm{cap}}(A,B)$ is the conformal capacity of the condenser $(A,B)$ (see [D1],[V]).
From (4.1) it follows (cf.~[D1, Sec.~4.6]) that there exists a constant $\rho'(\ell,\rho)$, depending only on $\ell$ and $\rho$, such that if $\phi$ is a Blaschke product, $\phi(0)=0$, $\text{\rm{val}}(\phi)\le \ell$, and the images of the branch points of $\phi$ are in $\{|z|\le \rho\}$, then the branch points of $\phi$ are in $\{|z|\le \rho' \}$.
Recall that the function $\bar f$ has a decomposition $\bar f = \phi\circ h$, where $\phi$ is a Blaschke product, $\text{\rm{val}}(\phi)\le \ell$, and $h$ is a self-homeomorphism of $\bold D$.
Since $S_{\tilde \p P}(v_0) = \smap P(v_0)=0$, we can assume that $\phi(0)=0$ and $h(0)=0$.
Thus $S_{\tilde \p P}(\brv P) \subset \{|z|\le \rho' \}$ because $\smap P(\brv P)\subset \{|z|\le \rho\}$, that is the images of the branch points of $\phi$ are in $\{|z|\le \rho\}$.
Hence, by (2), $\bar f$ is $K$-quasiregular, $K=K(d,\ell,\bar \rho)$, where $\bar \rho = \max\{\rho,\rho'(\ell,\rho)\}$.
From the definition of $\bar f$, (4.1), and the fact that $\bar f(\bold D)=\bold D$ we have 
$$
\text{\rm{cap}}\bigl( \bold D, S_{\tilde \p P}(\text{\rm{st}}(v_0)) \bigr) \le 
K\, \ell \, \text{\rm{cap}}\bigl( \bar f(\bold D), \bar f(S_{\tilde \p P}(\text{\rm{st}}(v_0))) \bigr)=
K \, \ell \, \text{\rm{cap}}\bigl( \bold D, S_{\p P}(\text{\rm{st}}(v_0)) \bigr),
$$
where $\text{\rm{st}}(v)$ is the union of all faces in \mcompk\ having $v$ as the common vertex.
The left-hand side inequality in (3) now follows from (1) and the fact that $C_{\tilde \p P}(v)\subset S_{\tilde \p P}(\text{\rm{st}}(v))$.
\qed
\enddemo

\remark{Remark~4.5}
Notice that the left-hand side inequality in (3) is the inverse result to the Discrete Schwarz Lemma of [BSt2] and [DSt].
\endremark

\heading 5.~Proof of Lemma~3.1 \endheading

We will now show Lemma~3.1.

\demo{Proof of Lemma~3.1}
Suppose that the assertion of the lemma is not true.
Then there exist a sequence of \blp s $\{\vp Pn\}$ satisfying 1)~--~3) and a sequence of pairs $\{v_n,w_n\}$ such that $v_nw_n$ is an edge in \mvcompk n and $r_{\vp Pn}(v_n)/r_{\vp Pn}(w_n)\ge n$, where \mvcompk n is the complex of \mvp Pn.
Let $\tilde \vp Pn$ be a univalent \blp\ for \mvcompk n.
Let $\bar f_n$ be the extension of the \dbp\ $f_n$ from $\tilde \vp Pn$ to \mvp Pn.
For $v\in \bdv n$ we define 
$$
A_n(v):= D_{\tilde \vp Pn}(v) \cup \bigcup \Sb \triangle \in \vsimpk n2\\ v\in \triangle \endSb S_{\tilde \vp Pn}(\triangle) \bigcup \Sb u \in \vsimpk n0\\ u\sim v \endSb A_{\tilde \vp Pn}(u,v),
$$
where the notation is as in Section~2.
From the definition of $\bar f_n$ it follows that $\bar f_n|_{A_n(v)}$ is 1-to-1.
Since $\bar f_n:\bold D \to \bold D$ is a branched covering, $\bar f_n(A_n(v))$ cannot contain images of all branch points of $\bar f_n$, i.e. $\bar f_n(\text{\rm{br}}_V(\vp Pn)) \setminus \bar f_n(A_n(v)) \ne \emptyset$.
In particular, $\bar f_n(\text{\rm{br}}_V(\vp Pn)) \setminus C_{\vp Pn}(v) \ne \emptyset$ because $C_{\vp Pn}(v)\subset \bar f_n(A_n(v))$.
This and the condition 3) imply
$$
r_{\vp Pn}(v)\le \tfrac12 (1+\rho) \quad \text{for $v\in \bdv n$}.
\tag{5.1}
$$

\remark{{\bf Note}} In the remainder of the proof the existence and convergence of various quantities will be taken for granted; this is due to the fact that one can always apply diagonalization techniques or take a subsequence if necessary because $\text{\rm{deg}}(\vcompk n)\le d$ for all $n$.
\endremark

We may assume without loss of generality (see Note) that there exists $m$ such that for each $n$ there are vertices $w^1_n,\dots , w^m_n, u_n$ with the following properties:
\roster
\item"{P(1)}"
$w^1_n:=w_n$,
\item"{P(2)}"
$w^1_n,\dots , w^m_n, u_n$ are consecutive neighbors of $v_n$ in \mvcompk n listed in the clockwise order,
\item"{P(3)}"
$\lim_{n\to \infty}\tfrac{r_n(v_n)}{r_n(w^i_n)} = \infty$, $i=1,\dots ,m$, and $\lim_{n\to \infty}\tfrac{r_n(v_n)}{r_n(u_n)} < \infty$, where $r_n(\cdot):=r_{\vp Pn}(\cdot)$.
\endroster
The following two observations can be helpful in seeing that the property~P(3) can be met.
\roster
\item"{(a)}"
If $v_n\in \bdv n$ for every $n$ then, since $\text{\rm{deg}}(\vcompk n)\le d$ and (5.1), 
$$
\varlimsup_{n\to \infty}\max \Sb v\in \vsimpk n0\\ v\sim v_n \endSb \{r_n(v)/r_n(v_n)\} >0,
$$ 
otherwise the circles $\{ C_{\vp Pn}(v) \}_{v:\, v\sim v_n}$ would not chain around $C_{\vp Pn}(v_n)$ from $\partial \bold D$ to $\partial \bold D$. (Figure~3(a).)
Moreover, $w_n$ cannot be the last of the consecutive neighbors of $v_n$ in \mvsimpk n0 listed in the counterclockwise order, otherwise there would exist $m\ge 1$ such that $\tilde w^1_n:= w_n, \tilde w^2_n, \dots , \tilde w^m_n, \tilde u_n$ are consecutive neighbors of $v_n$ in \mvsimpk n0 listed in the counterclockwise order, $\lim_{n\to \infty}r_n(v_n)/r_n(\tilde w^i_n) = \infty$, $i=1,\dots ,m$, and $\lim_{n\to \infty}r_n(v_n)/r_n(\tilde u_n) < \infty$; this would imply, however, that $C_{\vp Pn}(\tilde u_n)$ must intersect $\{|z|>1\}$ for sufficiently large $n$.
This is impossible.
(Figure~3(b).)
\item"{(b)}"
Similarly, if $v_n\in \interv n$ for every $n$ then, because $\text{\rm{deg}}(\vcompk n)\le d$, 
$$
\varlimsup_{n\to \infty}\max \Sb v\in \vsimpk n0\\ v\sim v_n \endSb \{r_n(v)/r_n(v_n)\} >0,
$$ 
otherwise the circles $\{ C_{\vp Pn}(v) \}_{v:\, v\sim v_n}$ would not wrap around $C_{\vp Pn}(v_n)$.
\endroster

\midinsert
\hskip13.0truecm $C_{\vp Pn}(\tilde w^i_n)$
\vskip-0.4truecm\centerline{\epsfysize=6.5truecm\epsffile[117 207 495 585]{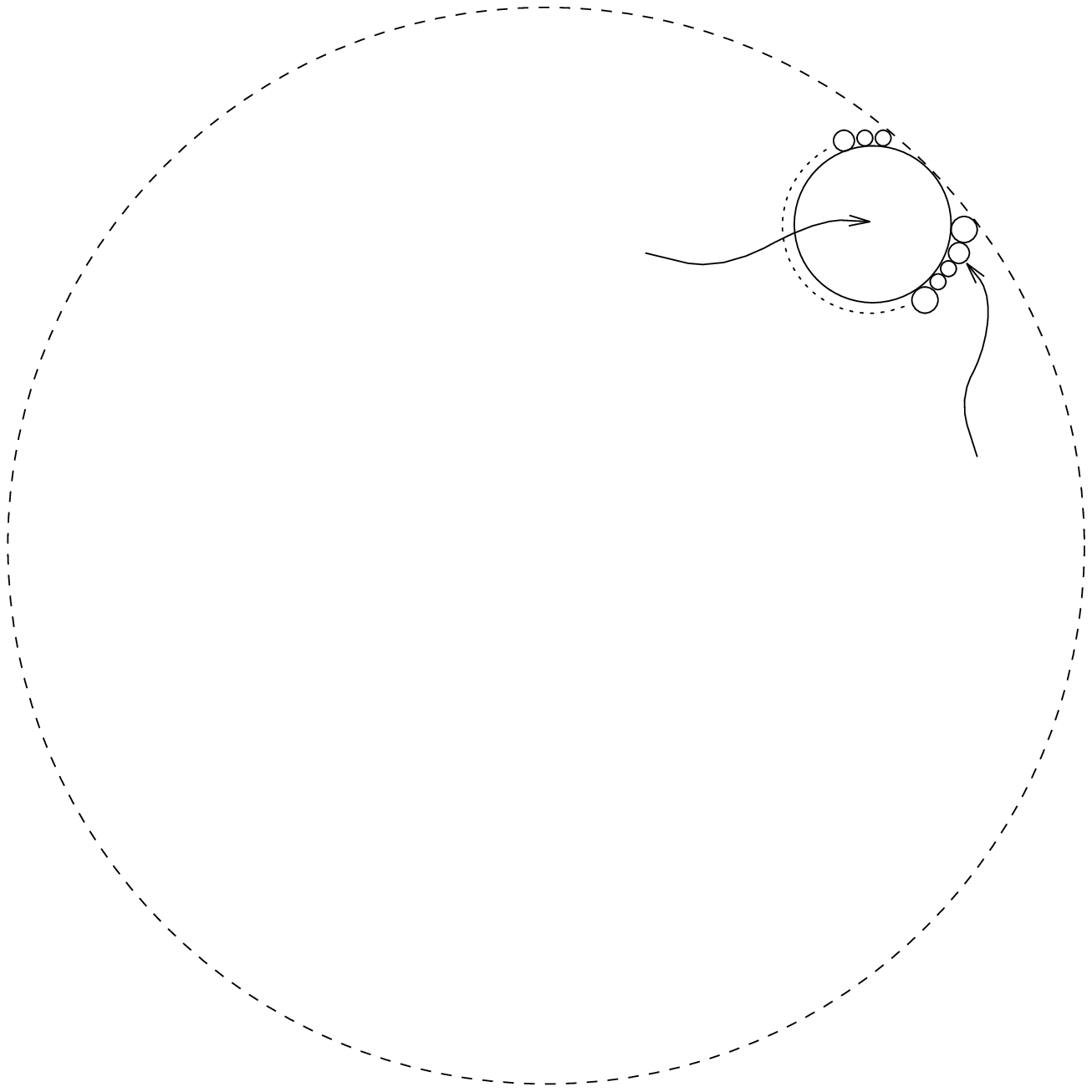} \hskip1.0truecm \epsfysize=6.5truecm\epsffile[117 207 495 585]{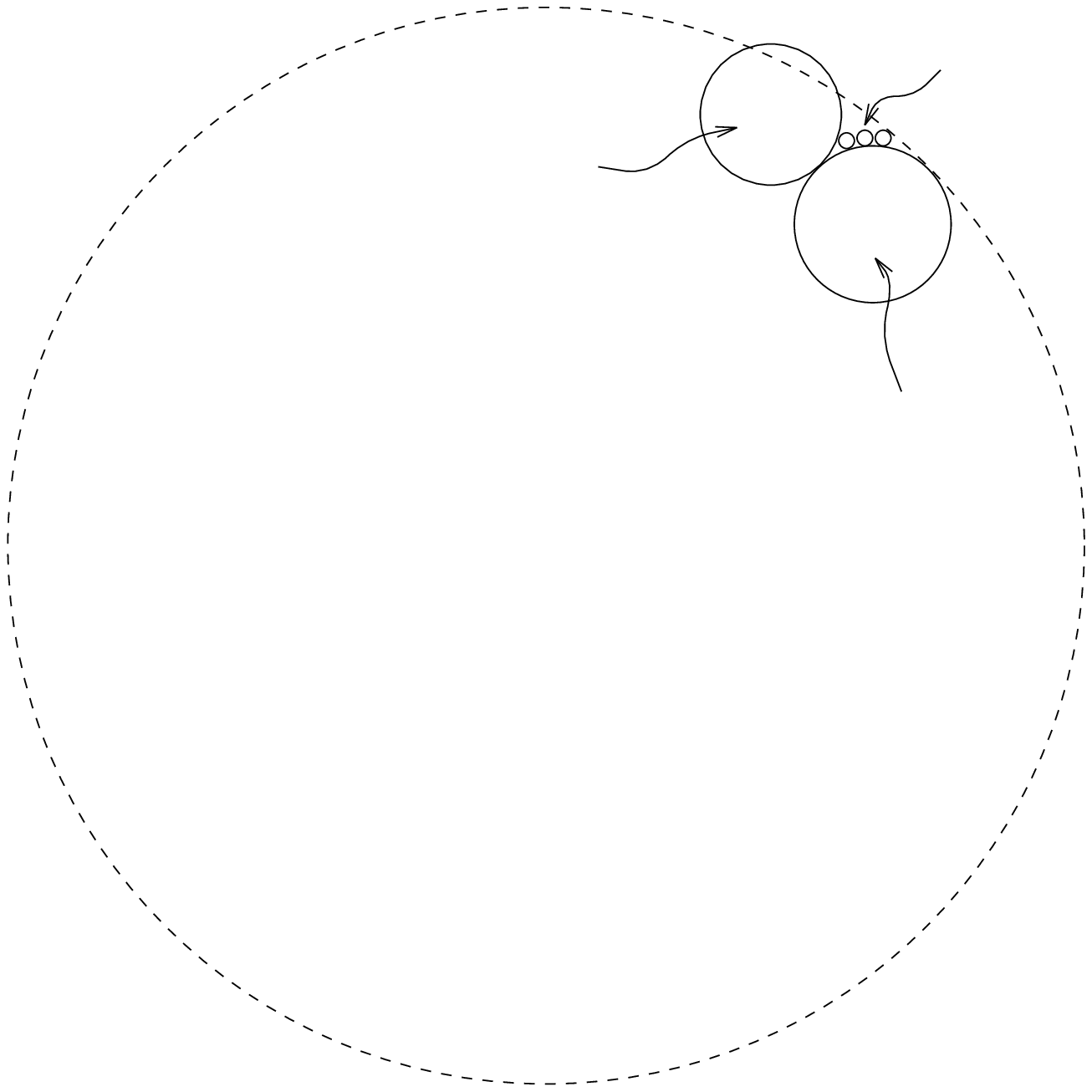}}
\vskip-5.9truecm\hskip9.5truecm $C_{\vp Pn}(\tilde u_n)$
\vskip0.0truecm\hskip2.35truecm $C_{\vp Pn}(v_n)$
\vskip0.6truecm\hskip12.0truecm $C_{\vp Pn}(v_n)$
\vskip-0.1truecm\hskip5.2truecm $C_{\vp Pn}(v)$
\vskip2.4truecm\hskip3.2truecm $\vp Pn$ \hskip7.0truecm $\vp Pn$
\vskip0.8truecm\hskip3.2truecm (a) \hskip7.0truecm (b)
\captionwidth{16pc}
\botcaption{Figure~3} The impossible.
\endcaption
\endinsert

We will now describe an inductive process.
We are going to define a double-sequence of sets of vertices $\{ w^1_n(j),\dots , w^{m(j)}_n(j), u_n(j) \}\subset \vsimpk n0$, $n=1,2,\dots$, and $j=1,2,\dots$, where $m(j)$ depends only on $j$.

Set $m(1):=m$, $w^i_n(1):= w^i_n$ for $i=1,\dots , m(1)$, and $u_n(1):= u_n$. For purely notational purposes we also set $u_n(0):= v_n$.
Suppose the sequences 
$$
\{ w^1_n(j),\dots , w^{m(j)}_n(j), u_n(j) \}_{n=1}^{\infty}
$$ 
have been defined for $j=1,\dots, N$.
We now define the sequence 
$$\{ w^1_n(N+1),\dots , w^{m(N+1)}_n(N+1), u_n(N+1) \}_{n=1}^{\infty}
$$
to be any sequence such that
\roster
\item"{$\wp(1)$}"
$w^1_n(N+1):=w^{m(N)}_n(N)$,
\item"{$\wp(2)$}"
$ w^1_n(N+1),\dots , w^{m(N+1)}_n(N+1), u_n(N+1)$ are consecutive neighbors of $u_n(N)$ in \mvcompk n listed in the clockwise order,
\item"{$\wp(3)$}"
$\lim_{n\to \infty}\tfrac{r_n(v_n)}{r_n(w^i_n(N+1))} = \infty$, $i=1,\dots ,m(N+1)$, and $\lim_{n\to \infty}\tfrac{r_n(v_n)}{r_n(u_n(N+1))} < \infty$.
\endroster
The existence of such a sequence easily follows from the fact that 
$$\lim_{n\to \infty}\tfrac{r_n(u_n(N))}{r_n(w^{m(N)}_n(N))} = 0, \ \lim_{n\to \infty}\tfrac{r_n(u_n(N-1))}{r_n(w^{m(N)}_n(N))} = 0,
$$
and the following slight modifications of the earlier observations (see also Note).
\roster
\item"{($\text{a}'$)}"
If $u_n(N)\in \bdv n$ for every $n$ then, since $\text{\rm{deg}}(\vcompk n)\le d$ and (5.1),
$$
\varlimsup_{n\to \infty}\max \Sb v\in \vsimpk n0\setminus \{u_n(N-1)\}\\ v\sim u_n(N) \endSb \{r_n(u_n(N))/r_n(v)\} >0.
\tag{5.2}
$$
For if the above were not true and $\lim_{n\to \infty}\tfrac{r_n(u_n(N-1))}{r_n(u_n(N))} > 0$ then, in order to preserve tangencies in \mvp Pn among the circles $C_{\vp Pn}(v)$, $v\sim u_n(N)$, and the circle $C_{\vp Pn}(u_n(N))$, $C_{\vp Pn}(u_n(N-1))$ would have to intersect $\{|z|>1\}$ for sufficiently large $n$, which is impossible.
When $\lim_{n\to \infty}\tfrac{r_n(u_n(N-1))}{r_n(u_n(N))} = 0$ the arguments of (a) can be used to show (5.2).
\item"{($\text{b}'$)}"
If $u_n(N)\in \interv n$ for every $n$ then it follows from the arguments of (b) that (5.2) holds in this case as well.
\endroster
It may a priori happen that after certain number of induction steps repetitions may start occuring.
We will show that given $L$ there exists $M_L$ and a subsequence of indices $\{n_k\}$ such that for each $k$ the number of different vertices (of $\compk_{n_k}$) in $\{ u_{n_k}(0),\dots, u_{n_k}(M_L) \}$ is at least $L$.

Suppose that this is not true.
Write $U_n(k,l)$ for the set $\{ u_n(k),\dots, u_n(l) \}$ and $\# U_n(k,l)$ for the number of different vertices in $U_n(k,l)$.
Then there exists $L$ such that for every $M$ there is $n_M$ such that $\# U_n(k,l)\le L$ for all $n\ge n_M$.
In particular, $U_n(0,M)\subseteq B_n(v_n,L)$ for all $n\ge n_M$, where $B_n(v,\varrho)$ denotes the ball in \mvsimpk n0 with center $v$ and radius $\varrho$ for the combinatorial metric in \mvsimpk n0.
Moreover, we can assume (see Note) that the double-sequence $\{u_n(j)\}_{j,n=1}^{\infty}$ has the following properties:
\roster
\item
there are indices $j_i$, $l_i$, and $\eta_i$, $i=1,2,...$, satisfying 
$$
\gather
0\le j_1< j_2< \dots, \quad 1\le l_1, l_2,\dots \le L-1, \quad j_i+l_i<j_{i+1}+l_{i+1},\\
\eta_1\le \eta_2\le\dots, \quad \text{and}\quad \eta_i\ge n_i \ \text{($n_i$ as above)},
\endgather
$$
such that for every $n\ge \eta_i$, $\# U_n(j_i,j_i+l_i) = l_i+1$ and $u_n(j_i)=u_n(j_i+l_i+1)$ (i.e., $\Gamma_n(i):= u_n(j_i)u_n(j_i+1)\cup \dots \cup u_n(j_i+l_i)u_n(j_i+l_i+1)$ is a simple closed edge-path in \mvcompk n of length at most $L$), and $\# U_n(j_i+1,j_{i+1}+l_{i+1}) =j_{i+1} + l_{i+1} - j_i$ (i.e., $u_n(j_i+1)u_n(j_i+2)\cup \dots \cup u_n(j_{i+1}+l_{i+1}-1)u_n(j_{i+1}+l_{i+1})$ is a simple edge-path in \mvcompk n)
\item
for all $n\ge \eta_i$ the number of vertices in $B_n(u_n(0),L)\cap \vsimpk n0$ enclosed by $\Gamma_n(i)$ is constant, say $\chi(i)$.
\endroster
Notice that (2) follows from the fact that the number of vertices of \mvcompk n contained in $B_n(v_n,L)$ (recall, $u_n(0)=v_n$) is at most $Ld^L$ because $\text{\rm{deg}}(\vcompk n)\le d$ for all $n$.
Also, observe that ($\text{a}'$) and ($\text{b}'$) imply that each $l_i$ must be at least 2.

Let $\iota:= \min\bigl\{ i: \chi(i)=\min_j\{\chi(j)\}\bigr\}$.
Write $\Gamma_n$ for $\Gamma_n(\iota)$.
Let \mvcomp On be the subcomplex of \mvcompk n bounded by $\Gamma_n$.
(The path $\Gamma_n$ is included in \mvcomp On.)
From $\wp(1)$, $\wp(2)$, and the induction definition it follows that for each $n$ one of the following two cases occurs.
\roster
\item"{Case~1:}"
$$
\biggl( \bigcup_{j=j_{\iota}+1}^{l_{\iota}}\{ w^i_n(j) \}_{i=1}^{m(j)} \biggr) \cap \vcomp On =\emptyset;
$$
\item"{Case~2:}"
$$
\bigcup_{j=j_{\iota}+1}^{l_{\iota}}\{ w^i_n(j) \}_{i=1}^{m(j)} \subseteq \vcomp On.
$$
\endroster
By taking a subsequence if necessary, we can assume that either Case~1 or Case~2, but not both, holds for all $n$ (see Fig.~4).

\midinsert
\centerline{\epsfysize=6.5truecm\epsffile[75 260 533 530]{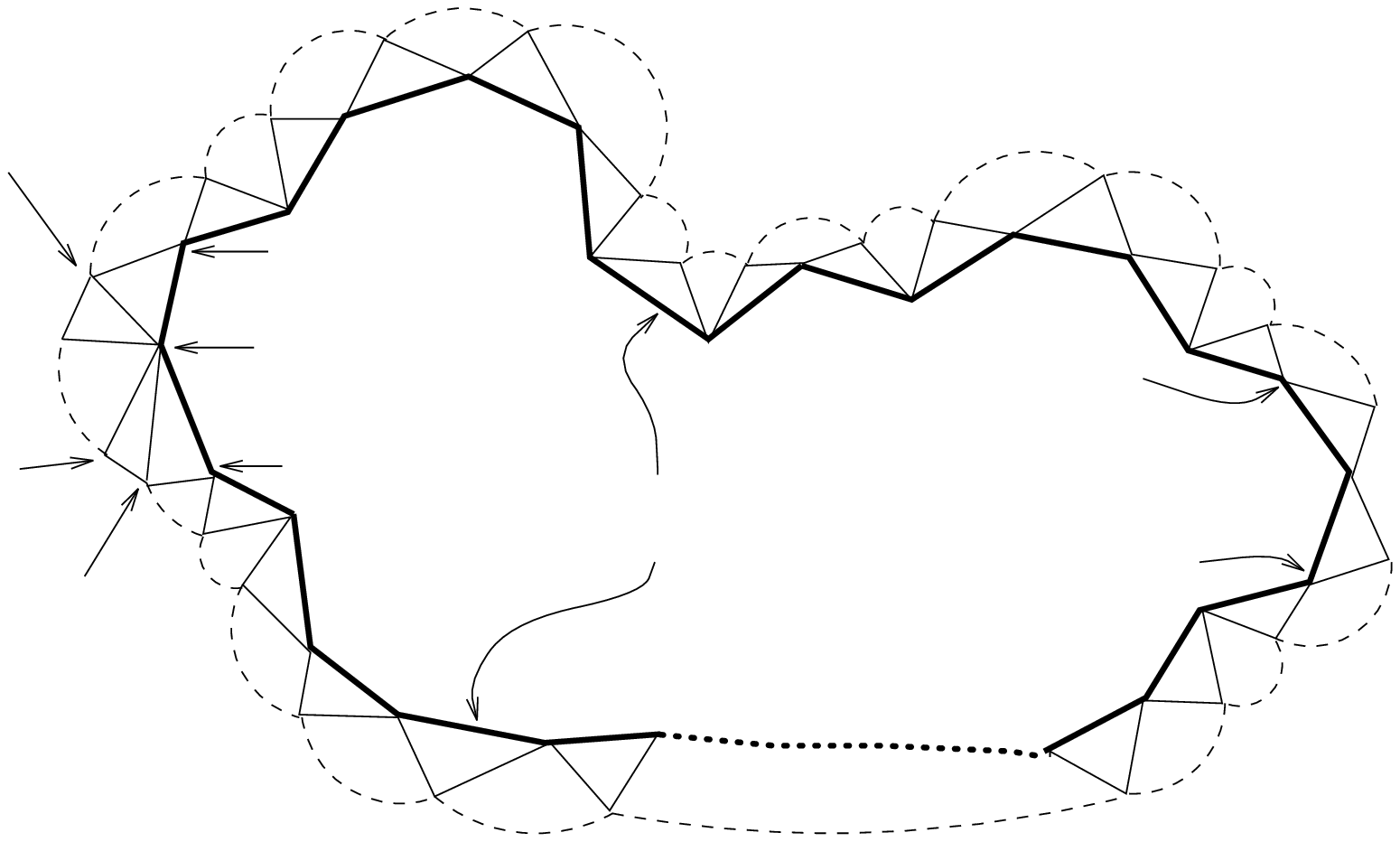}}
\vskip-5.8truecm\hskip1.0truecm $w^{m(j)}_n(j)$
\vskip0.3truecm\hskip4.0truecm $u_n(j+1)$
\vskip1.2truecm\hskip0.85truecm $w^2_n(j)$
\vskip-1.5truecm\hskip3.9truecm $u_n(j)$
\vskip-0.3truecm\hskip8.65truecm $u_n(j_{\iota}+l_{\iota})$
\vskip0.2truecm\hskip4.1truecm $u_n(j-1)$
\vskip-0.05truecm\hskip6.7truecm $\Gamma_n$
\vskip-0.9truecm\hskip10.9truecm $u_n(j_{\iota})$ \hskip0.2truecm $u_n(j_{\iota}+l_{\iota}+1)$
\vskip0.7truecm\hskip1.6truecm $w^1_n(j)$
\vskip-1.0truecm\hskip9.15truecm $u_n(j_{\iota}+1)$
\vskip2.0truecm\centerline{Case~1}
\vskip0.0truecm\centerline{\epsfysize=7.0truecm\epsffile[75 260 533 530]{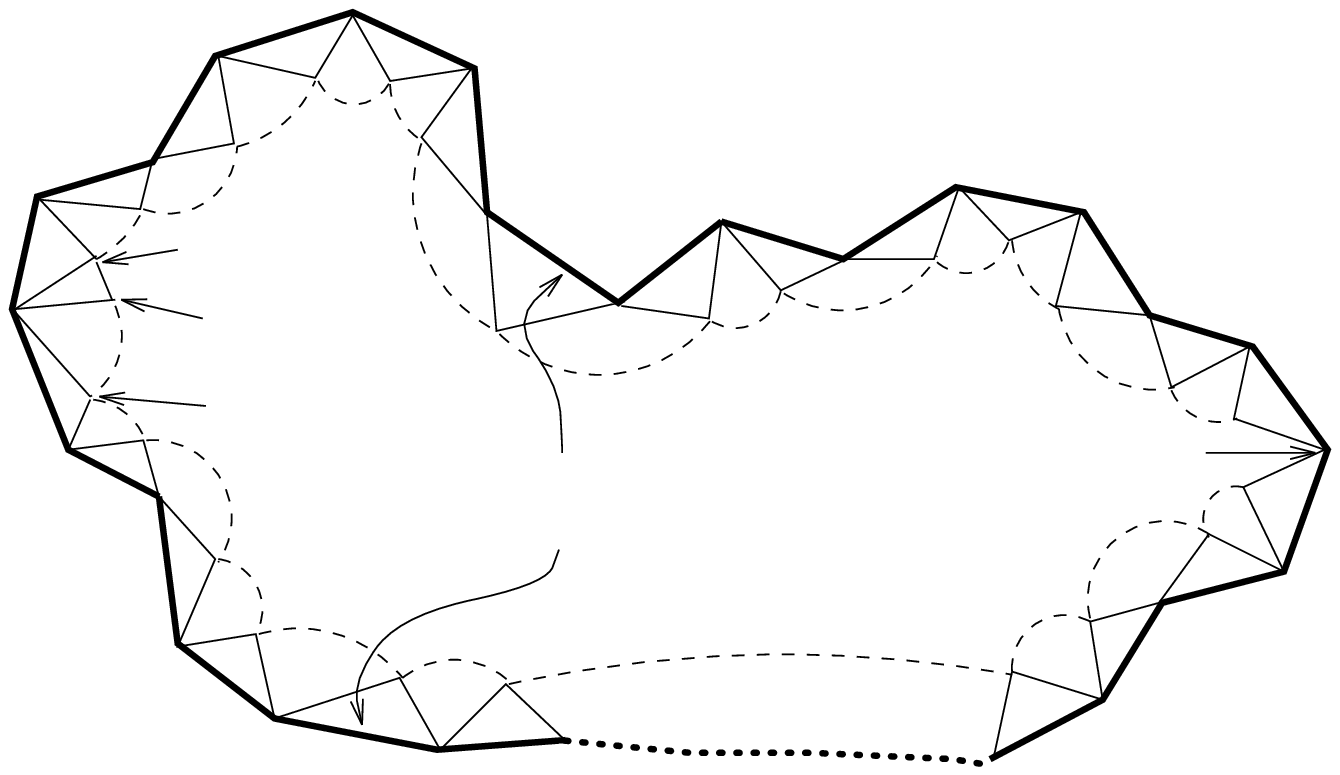}}
\vskip-5.5truecm\hskip0.95truecm $u_n(j-1)$
\vskip-0.05truecm\hskip3.85truecm $w^1_n(j)$
\vskip-0.08truecm\hskip1.5truecm $u_n(j)$
\vskip-0.35truecm\hskip4.05truecm $w^2_n(j)$
\vskip-0.47truecm\hskip11.9truecm $u_n(j_{\iota}+1)$
\vskip0.0truecm\hskip4.05truecm $w^{m(j)}_n(j)$
\vskip0.23truecm\hskip6.5truecm $\Gamma_n$
\vskip-0.9truecm\hskip10.3truecm $u_n(j_{\iota})$ \hskip0.9truecm $u_n(j_{\iota}+l_{\iota}+1)$
\vskip-0.4truecm\hskip1.1truecm $u_n(j+1)$
\vskip0.35truecm\hskip12.15truecm $u_n(j_{\iota}+l_{\iota})$
\vskip1.5truecm\centerline{Case~2}
\captionwidth{16pc}
\botcaption{Figure~4} The cases.
\endcaption
\endinsert

We now recall (cf. [D1]) that if \mcompk\ is a finite triangulation of a disc then the following equation, which is the Euclidean version of Gauss-Bonnet formula, describes the relationship between interior and boundary angle sums induced by a \cp\ \mp P for \mcompk:
$$
\sum_{v\in \intk} (2\pi - \Theta_{\p P}(v)) = 2\pi - \sum_{v\in \bdk} (\pi - \gamma_{\p P}(v)).
\tag{5.3}
$$

Let
$$
\multline
\beta_n(j):= \alpha_{\p P_n}\bigl( u_n(j), \langle u_n(j),u_n(j-1),w^1_n(j) \rangle \bigr) +
\\
\sum_{i=1}^{m(j)-1} \alpha_{\p P_n}\bigl( u_n(j), \langle u_n(j),w^i_n(j),w^{i+1}_n(j) \rangle \bigr) +
\\ 
\alpha_{\p P_n}\bigl( u_n(j), \langle u_n(j),u_n(j+1),w^{m(j)}_n(j) \rangle \bigr),
\endmultline
$$
and $\gamma_n(j):= \sum_{\triangle \in \vcomp On} \alpha_{\p P_n}(u_n(j), \triangle)$, i.e. $\gamma_n(j)$ is the boundary angle sum at $u_n(j)$ induced by $\p P_n|_{\vcomp On}$.

Suppose that Case~1 holds for all $n$.
Then $u_n(j_{\iota}+1),\dots, u_n(j_{\iota}+l_{\iota}) \in \interv n$.
Hence $\gamma_n(j)\ge 2\pi - \beta_n(j)$ for $j=j_{\iota}+1,\dots, j_{\iota}+l_{\iota}$ and all $n$.
By (5.3),
$$
0\ge \sum_{v\in \text{\rm{int}}\, \vcomp On^0} (2\pi - \Theta_{\p P_n}(v)) =
2\pi - \sum_{j=j_{\iota}}^{l_{\iota}} (\pi - \gamma_n(j)) \ge
\pi + \gamma_n(j_{\iota}) + \sum_{j=j_{\iota}+1}^{l_{\iota}} (\pi - \beta_n(j)).
\tag{5.4}
$$
From $\wp(3)$ we have
$$
\lim_{n\to \infty} \frac{r_n(u_n(j))}{r_n(w^i_n(j))} =\infty \quad \text{for $i=1,\dots , m(j)$ and any $j$}.
\tag{5.5}
$$
Hence
$$
\lim_{n\to \infty} \beta_n(j) = 0 \quad \text{for every $j$}.
\tag{5.6}
$$
This and (5.4) give a contradiction.
Thus Case~1 cannot occur.

Suppose Case~2 holds.
Then $\wp(1)$ and $\wp(2)$ imply that either $u_n(j_{\iota}+l_{\iota}+2) \in \text{\rm{int}}\, \vcomp On^0$ or $u_n(j_{\iota}+l_{\iota}+2) = u_n(j_{\iota}+1)$.
By choosing a subsequence if necessary, we may assume that one of these two possibilities holds for all $n$.
We will show that the first possibility must be ruled out.

If $u_n(j_{\iota}+l_{\iota}+2) \in \text{\rm{int}}\, \vcomp On^0$ for all $n$ then it follows from $\wp(1)$ -- $\wp(3)$ that the path $u_n(j_{\iota}+l_{\iota}+1)u_n(j_{\iota}+l_{\iota}+2) \cup\dots\cup u_n(j_{\iota+1}+l_{\iota+1})u_n(j_{\iota+1}+l_{\iota+1}+1)$ must stay in $\text{\rm{int}}\, \vcomp On^0$ for every $n\ge \eta_{\iota+1}$.
(Notice that for all sufficiently large $n$ (see Note), say $n\ge \eta_{\iota+1}$, 
$$
\{u_n(j_{\iota}+l_{\iota}+1),\dots, u_n(j_{\iota+1}+l_{\iota+1}+1)\} \cap \biggl( \bigcup_{j=j_{\iota}+1}^{j_{\iota}+l_{\iota}+1}\{ w^i_n(j) \}_{i=1}^{m(j)} \biggr)  = \emptyset
$$
because of $\wp(1)$ -- $\wp(3)$, and either $u_n(j_{\iota+1}+l_{\iota+1}+1) \in \text{\rm{int}}\, \vcomp On^0$, in which case $\Gamma_n(\iota+1)$ is contained in the interior of \mvcomp On, or $u_n(j_{\iota+1}+l_{\iota+1}+1) = u_n(j_{\iota}+l_{\iota}+1)$, in which case $j_{\iota+1} = j_{\iota}+l_{\iota}+1$ and $\Gamma_n(\iota+1)$ is the path of the type considered in Case~1, i.e. cannot exist.)
Now, since $\Gamma_n(\iota+1)\subset B_n(u_n(0),L)$ for every $n$, it is easy to see that $\chi(\iota+1) < \chi(\iota)$, and we get a contradiction with the definition of $\iota$.

Thus we are left with the second possibility, i.e. $u_n(j_{\iota}+l_{\iota}+2) = u_n(j_{\iota}+1)$.
In this case $j_{\iota+1} = j_{\iota}+1$ and $l_{\iota+1} = l_{\iota}$, and (5.3) implies
$$
\sum_{v\in \text{\rm{int}}\, \vcomp On^0} (2\pi - \Theta_{\p P_n}(v)) =
2\pi - \sum_{j=j_{\iota+1}+1}^{j_{\iota+1}+l_{\iota+1}+1} (\pi - \beta_n(j)).
\tag{5.7}
$$
Since $\p P_n$'s are packings, by Theorem~4.1,
$$
l_{\iota+1}+1 \ge 3 + 2\left( \sum_{v\in \text{\rm{br}}_V(\p P_n)\cap \text{\rm{int}}\, \vcomp On^0} \text{\rm{ord}}_{\p P_n}(v) \right).
$$
Hence, (5.7) and (5.5) imply the following contradiction
$$
\multline
-2\pi\left( \sum_{v\in \text{\rm{br}}_V(\p P_n)\cap \text{\rm{int}}\, \vcomp On^0} \text{\rm{ord}}_{\p P_n}(v) \right) = \sum_{j=j_{\iota+1}+1}^{j_{\iota+1}+l_{\iota+1}+1} (\pi - \beta_n(j)) \underset{n\to\infty}\to\longrightarrow 2\pi - (l_{\iota+1}+1)\pi \le \\
-\pi -2\pi\left( \sum_{v\in \text{\rm{br}}_V(\p P_n)\cap \text{\rm{int}}\, \vcomp On^0} \text{\rm{ord}}_{\p P_n}(v) \right).
\endmultline
$$
Thus we have shown that Case~2 cannot occur neither.
Therefore for any given $L$ there exists $M_L$ and a subsequence of indices $\{n_k\}$ such that for each $k$ the number of different vertices (of $\compk_{n_k}$) in $\{u_{n_k}(0),\dots, u_{n_k}(M)\}$ is at least $L$.

Let $L:=2(\ell +1)$.
We may assume that $\#U_n(0,M_L)\ge L$ for every $n$.
We will show that for each sufficiently large $n$ there is a point $x_n\in \bold D$ which is covered by at least $\ell +1$ different discs from $\{D_n(0),\dots, D_n(M_l)\}$, where $D_n(j):= D_{\p P_n}(u_n(j))$.
This, however, will imply that packings $\p P_n$ are of valence no less then $\ell +1$ in contradiction to our assumption~2), and the proof of the lemma will be complete.

We now show the existence of $x_n$'s.
First we ``multiply'' $\bold D$ (i.e., scale and rotate) by the factor $\zeta_n\in\bold C$ so that $\zeta_nC_n(0)$ has the radius equal to 1 and so that the tangency point of $\zeta_nC_n(0)$ and $\zeta_nC_n(1)$ is the point $1\in \bold C$, where $C_n(j):= C_{\p P_n}(u_n(j))$.
From $\wp(3)$ it follows that there exists $\delta>0$ such that for each $j$, $j=1,\dots, M_L$, the radius of $\zeta_nC_n(j)$ is at least $\delta$ for all sufficiently large $n$.
Moreover, $\wp(1)$ -- $\wp(3)$ imply that $$
\bigl\{ \zeta_nC_n(w^1_n(1)),\dots,\zeta_nC_n(w^{m(1)}_n(1)),\cdots , \zeta_nC_n(w^1_n(M_L)),\dots,\zeta_nC_n(w^{m(M_L)}_n(M_L)) \bigr\}
$$
is a chain of circles such that, as $n\to\infty$, it collapses geometrically to a point in $\bold C$ which is 1.
From the construction of the chain it follows that the tangency points of pairs $\{\zeta_nC_n(j),\zeta_nC_n(j+1)\}$, $j=0,\dots, M_L-1$, have the common limit point in $\bold C$, as $n\to\infty$, which is 1.
As each $\zeta_nC_n(j)$, $j=0,\dots, M_L-1$, has its radius no smaller then $\delta$, it is easy to see that for each sufficiently large $n$ there is a point $y_n\in\bold C$ which is covered by at least a half of different discs from $\{\zeta_nD_n(0),\dots, \zeta_nD_n(M_l)\}$, i.e. by at least $L/2=\ell +1$ of different discs.
We now take $x_n:=y_n/\zeta_n$ and the proof is finished.
\qed
\enddemo

\enddocument